\documentclass[11pt]{article}
\usepackage{amsmath,amssymb,amsbsy,amsfonts,amsthm,latexsym,
                          amsopn,amstext,amsxtra,euscript,amscd}
\newtheorem{thm}{Theorem}
\newtheorem{theorem}[thm]{Theorem}

\newtheorem{corollary}[thm]{Corollary}
\newtheorem{conjecture}[thm]{Conjecture}

\def\\{\cr}
\def\({\left(}
\def\){\right)}
\def\[{\left[}
\def\]{\right]}
\def\<{\langle}
\def\>{\rangle}

\def\cD{\mathcal D}
\def\cE{\mathcal E}

\def\cL{\mathcal L}

\def\cP{\mathcal P}

\def\cS{\mathcal S}

\def\Z{\mathbb{Z}}

\begin{document}

\title{\bf Complexity of  Inverting the Euler Function}

\author{
{\sc Scott Contini} \\
{Department of Computing}\\
{Macquarie University} \\
{Sydney, NSW 2109, Australia} \\
{\tt contini@ics.mq.edu.au}
\and
{\sc Ernie Croot} \\
{ School of Mathematics} \\
{ Georgia Institute of Technology } \\
{ Atlanta, GA 30332, USA} \\
{\tt  ecroot@math.gatech.edu  } \\
\and
{\sc Igor E.~Shparlinski} \\
{Department of Computing}\\
{Macquarie University} \\
{Sydney, NSW 2109, Australia} \\
{\tt igor@ics.mq.edu.au}
}
\date{\today}

\date{\today}

\pagenumbering{arabic}

\date{\today}
\maketitle

\begin{abstract}
We present an algorithm to invert the Euler function $\varphi(m)$.
The algorithm, for a given integer $n\ge 1$, in polynomial time ``on average'',
   finds the
set $\Psi(n)$ of all solutions $m$ to the equation $\varphi(m) =n$.
In fact, in the worst
case the set $\Psi(n)$ is exponentially large and cannot be
constructed by a polynomial
time algorithm. In the opposite direction, we show, under some widely
accepted number theoretic conjecture, that the Partition Problem, an
{\bf NP}-complete problem, can be reduced, in polynomial time, 
to the problem of deciding whether 
$\varphi(m) = n$ has a solution, for polynomially (in the input size of
the Partition problem) many values of $n$.  In fact, the following
problem is {\bf NP}-complete:  Given a set of positive integers $S$, decide
whether there is an $n \in S$ satisfying $\varphi(m) = n$, for some integer
$m$.  Finally,  we establish close
links between the problem of inverting the Euler function  and the
integer factorisation problem.
\end{abstract}

\section{Introduction}

In this paper we study the complexity of a new number theoretic problem,
namely the complexity of  inverting the   {\it Euler function\/}  $\varphi(m)$,
which, as usual, for an integer $m \ge 1$, is defined by
$$
\varphi(m)=\#(\Z/m\Z)^\times=\prod_{p^\alpha\,\|\,m}p^{\alpha-1}(p-1).
$$

It is  widely believed that
computing the Euler function is equivalent to the integer
factorisation problem. Moreover, let  $\cP_2$ denote the set of
    positive integers $n$  which are
products of two distinct primes $p$ and $q$  (with the additional condition
    $p \equiv q \equiv 3 \pmod 4$ such numbers are often called {\it
Blum integers\/}).
Then for $n   \in \cP_2$ finding  $\varphi(n) $ is indeed equivalent
to factoring $n$.
Here we concentrate on the dual question of inverting the Euler
function, which apparently has not yet
been addressed in the literature.
More precisely,  given an integer $n\ge 1$, we  want to find the set
$\Psi(n)$ of all integer solutions $m\ge 1$ to the equation
$\varphi(m) = n$.

    Here we design
an algorithm which solves this problem in exponential time in $\log n$
in the worst case, and in polynomial time for ``almost all''
$n$ (provided the prime number factorisation of $n$ is given).
Because for infinitely many $n$ the cardinality
of $\Psi(n)$  is exponentially large, any algorithm for inverting $\varphi$
{\it must} run in time exponential in $\log n$ in the worst case (or nearly
exponential).
Indeed, from the proof of Theorem~4.6 of~\cite{Pom2} we see that
for infinitely many $n$,
$$
\# \Psi(n) \ge  n^{\gamma +o(1)}
$$
where $\gamma> 0$ is any constant such that for any sufficiently
large $X$ there are at least $X^{1 + o(1)}$ primes $p\le X$ such that all
prime divisors of $p-1$ are less than $X^{1-\gamma}$ (see also~\cite{Pom1}).
   By Theorem~1 of~\cite{BaHa}
one can take $\gamma = 0.7039$.

A natural question is whether the decision problem for inverting $\varphi$
is any easier.  We recall that $n$ is called a {\it totient\/},
if there exists an integer $m$ satisfying $\varphi(m) = n$.
Given an integer $n$, and its prime factorization,
how efficiently can we determine
whether $n$ is a totient?  Because the
output of any algorithm solving this problem need only be a single
bit, we cannot so easily say that the running time must be exponential in
$\log n$, as we did in the case of determining all the solutions $m$.
We  prove in Section~\ref{np_section}
the somewhat surprising result that, assuming a certain strong form
of the famous {\it Hardy--Littlewood\ } prime $k$-tuplet conjecture
(in the case $k=2$), there is a polynomial time reduction
from the Parition Problem, an {\bf NP}-complete problem, to the 
question of whether $\varphi(m) = n$ has a solution for a certain small
set of integers $n$.  In particular this shows that the following
problem is {\bf NP}-complete (assuming the Hardy--Littlewood conjecture):
Given a set of integers $\cS$, determine whether it contains a
totient.  Although at the present time the Hardy--Littlewood 
conjecture is out of reach, there are a number of results in this direction
which leave little doubt that the conjecture is correct, for example,
see~\cite{Bal}.

Furthermore, in Section~\ref{sec:Factor} we obtain an unconditional
reduction from the problem of factoring integers $n \in \cP_2$
to that of inverting the Euler function.
As we have remarked, any polynomial time
algorithm to compute the
Euler function leads to a factorization algorithm for
integers  of the form $n = pq$ where  $p$ and $q$ are  primes.
Here we prove  a somewhat dual statement by showing that
any polynomial time
algorithm to invert the Euler function (in the sense that the
running time is $(\# \Psi(n) + \tau(n) + \log n)^{O(1)} $ ),
where the factorization of
$n$ is {\it not} given, leads to a probabilistic
polynomial time factorisation algorithm for $n \in \cP_2$.
This result is certainly weaker than that of Section~\ref{np_section}
but is not based on any unproven assumptions.

The growth, distribution in arithmetic progressions
and in other special sets of elements of the values of the Euler function,
and many other similar questions,
have extensively been studied  in the literature,
see~\cite{BFPS,BLFS,DP,EP,Ford,FKP,Pom1,Pom2} and references therein.
Nevertheless the considered here questions seem to be new and have never
been  studied.  We also remark that analogues of our results can be obtained
for the sum of divisors function $\sigma(m)$ and for several more similar
number theoretic functions.

\section{Notation}

We use $\omega(m)$ and   $\tau(m)$  to denote the total number of
distinct prime and positive integer  divisors of a positive integer $m$,
respectively (we also define $\omega(1) = 0$, $\tau(1) =1$).

We also use the Vinogradov symbols
$\gg$, $\ll$, $\asymp$ as well as the Landau symbols $O$ and $o$
with their regular meanings  (we recall that
$U\ll V$ and  $U = O(V)$ are both equivalent to the
inequality $|U| \le cV$ with some constant $c> 0$ and $U \asymp V$ is
equivalent $U \ll V \ll U$).
The implied constants in the
symbols $O$, $\gg$, $\ll$
and $\asymp$ are always absolute unless indicated otherwise.

\section{Constructing $\Psi(n)$}

Our algorithm to find $\Psi(n)$
makes use of the   prime power factorization of $n$.
If we were to modify our algorithm to find $\Psi(n)$ where
the factorisation of $n$ is not given,  but is first found
by using a probabilistic factoring algorithm (see~\cite{CrPom}), then
for most integers $n$,
factoring would dominate the overall complexity of the algorithm.
In the worst case, however, where $\Psi(n)$ is ``large'',
the running time of the rest of the algorithm would dominate this factoring
step.  For our algorithm, we simply assume that
the  prime number factorisation of $n$
is given, which has the additional advantage of making our
algorithm deterministic
(while making factoring $n$ a part of the
algorithm would make it probabilistic).

\begin{theorem}
\label{thm:Algorithm}
There exists a deterministic algorithm which given
the  prime number factorisation
$$
n =  p_1^{\alpha_1} \ldots p_s^{\alpha_s}
$$
of an integer $n\ge 2$,  constructs $\Psi(n)$  in time
$$
T(n) \le \(\#\Psi^*(n) + \tau(n) + \log n\)^{O(1)},
$$
where
$$
\Psi^*(n)=\bigcup_{d | n} \Psi(d).
$$
\end{theorem}

\begin{proof}  Basically, we give an algorithm which efficiently
finds all representations of
$n$ of the following type:
$$
n=\prod_{j=1}^k \ell_j^{\gamma_j} (\ell_j - 1),
$$
where $\ell_1 < \cdots < \ell_k$ are primes, and
where $\gamma_1,\ldots,\gamma_k \geq 0$
are integers.  Each such representation corresponds to a solution
$\varphi(m) = n$, where
$$
m=\prod_{j=1}^k \ell_j^{\gamma_j + 1}.
$$

Our algorithm is iterative and builds a
graph, where all the vertices on the $j$th level correspond to a
certain list $\cE_j$, and
where each solution $m$ to $\varphi(m) = n$
corresponds to some path from a vertex back to the top list $\cE_1$,
although not all such paths  correspond to such a solution.
The vertices in each of these lists are  assigned a certain value,
which is  an ordered pair of the form $(\ell_j,\gamma_j)$, where
$\ell_j$ is prime, and where $\gamma_j \geq 0$ is an integer.

Given an integer $n$, we let $\cD(n)$ denote the set of
divisors of $n$.  If we are given the prime power factorisation of $n$,
then we can easily construct the set $\cD(n)$ in time $\tau(n)^{O(1)}$.

We now describe $\cE_1$:  We let $\cE_1$ be a set of vertices, one
for each ordered pair  $(\ell_1,\gamma_1)$,
where $\ell_1$ is a prime, and $\gamma_1 \geq 0$, such that the
number $e = \ell_1^{\gamma_1} (\ell_1-1)$ lies in $\cD(n)$; that is, $e | n$.
We also remark that for every $e$ there are at most two possible
pairs  $(\ell_1,\gamma_1)$.
The vertices in this list are not be linked to each other, but are
each   doubly
linked to entries of the yet to be mentioned list $\cE_2$.

The list $\cE_2$ is created as follows:  We scan through $\cE_1$, and for
each vertex $v$ in $\cE_1$, having the value
$(\ell_1, \gamma_1)$, we consider the integer
$n_0  = n/\ell_1^{\gamma_1}(\ell_1-1)$.  Then, among the integers
$d_0  \in \cD(n_0)$ (divisors of $n_0$), we locate all those corresponding
to vertices $v \in \cE_1$ having value $(\ell,\gamma)$, with
$\ell > \ell_1$; that is, $d_0 = \ell^\gamma (\ell -1)$.
In this way, we run through all the divisors $d$ of $n$
of the form
$$
d=\ell_1^{\gamma_1} (\ell_1 - 1) \ell^{\gamma} (\ell - 1),
\ \ \ell_1 < \ell\ \text{are prime}.
$$
The list $\cE_2$  then consists of one vertex for each of these different
ordered pairs $(\ell, \gamma)$, for each of the vertices $v \in \cE_1$;
and, this vertex is  doubly linked to its ancestor $v \in \cE_1$.

We note that each vertex in $\cE_2$ has a unique ancestor; and, different
vertices in $\cE_2$ may have the same value $(\ell_2,\gamma_2)$.

In general, suppose we have constructed the list $\cE_j$.  Then,
the list $\cE_{j+1}$ is constructed as follows:  By running through
the vertices $v \in \cE_j$, and
then considering the unique path from $v$ back to its ancestors in
$\cE_{j-1}, \cE_{j-2},\ldots,\cE_1$,
we get that these vertices (along with $v$)
correspond to a sequence of ordered pairs
$(\ell_j, \gamma_j),\ldots.,(\ell_1,\gamma_1)$, which represents a divisor $d$
of $n$ of the form
$$
d=\prod_{i=1}^j \ell_i^{\gamma_i} (\ell_i - 1),\ \ell_1 < \ell_2 < \cdots
< \ell_j\ {\rm are\ prime}.
$$
We let $n_0 = n/d$, and then we scan through the set
$\cD(n_0)$, looking for elements of the form $\ell^{\gamma}
(\ell - 1)$, where $\ell > \ell_j$ is a prime number.  We then
add a vertex to $\cE_{j+1}$, assign it the value
$(\ell, \gamma)$,
and doubly link it to the vertex $v \in \cE_j$.
After we have done this for all these ordered pairs
$(\ell,\gamma)$ generated by considering all $v \in \cE_j$,
the construction of $\cE_{j+1}$ is completed.

We continue constructing these lists, until we reach a list $\cE_t$
having no children.  Since $n$ has $O(\log n)$ prime power factors,
and since each new level in the graph corresponds to a string of
divisors $d_1,\ldots,d_t$ where $d_1\cdots d_t | n$, we conclude that
$t = O(\log n)$.

It is obvious that each path from a vertex back to $\cE_1$ along its unique
ancestors in the graph corresponds either to a proper divisor
$$
d=\prod_{j=1}^h \ell_j^{\gamma_j} (\ell_j - 1),
\ \ell_1 < \ell_2 < \cdots < \ell_h\ {\rm are\ prime}
$$
of $n$ such that there is no pair $(\ell, \gamma)$, $\ell > \ell_h$ prime,
$\gamma \geq 0$, with $\ell^\gamma (\ell - 1) d | n$; or, we have that
$d = n$.   Now, if $d = n$, then $d$ corresponds to the solution
$$
m=\prod_{j=1}^h  \ell_j^{\gamma_j+1}
$$
of $\varphi(m) = n$.  Thus, by considering the paths from vertices
corresponding to $n$ back to
$\cE_1$ corresponding to $d = n$, we obtain the set $\Psi(n)$.

Finally, it is obvious that the running time of the algorithm is
proportional to
$$
(L + \tau(n) + \log n)^{O(1)},
$$
where $L$ is the number of paths throughout the above graph
which is
$\# \Psi^*(n)$.
\end{proof}

To address the average performance of the algorithm, we require the
following bound:

\begin{theorem}
\label{thm: Psi^* aver}  The bound following bound holds:
$$
\sum_{n \leq x} \# \Psi^*(n) \ll  x \log x.
$$
\end{theorem}

\begin{proof}  We have that
\begin{equation} \label{psi_sum}
\sum_{n \leq x} \# \Psi^*(n) = \sum_{n \leq x} \sum_{d | n}
\# \Psi(d) \le  x \sum_{d \leq x} {\# \Psi(d) \over d}.
\end{equation}
Now,
$$
\sum_{d \leq x} \# \Psi(d)
=\#\{ n \geq 1\ :\ \varphi(n) \leq x\}
= \(\frac{\zeta(2)\zeta(3)}{\zeta(6)} + o(1)\)x,
$$
see~\cite{Bat}. So, by partial summation,
we conclude that
$$
\sum_{d \leq x} {\# \Psi(d) \over d}=O(\log x),
$$
which, together with~\eqref{psi_sum}  finishes the proof.
\end{proof}

An almost immediate corollary of  Theorem~\ref{thm: Psi^* aver}, together with
the well known bound
\begin{equation}
\label{eq:sums of tau}
\sum_{n \le x} \tau(n) = O(x\log x),
\end{equation}
see Theorem~2 in Section~I.3.2 of~\cite{Ten}, and
Theorem~\ref{thm:Algorithm}, is the following:

\begin{corollary}
   For every $A > 0$, there exists $B > 0$, so that
for all but at most $O(x/\log^A x)$ integers $n \leq x$ we have that
the algorithm in Theorem~\ref{thm:Algorithm}
finds $\Psi(n)$ in time $\log^B n$.
\end{corollary}

\section{{\bf NP}-completeness of Totient Testing}
\label{np_section}

A natural question is whether it is any easier to decide if, given $n$,
there exists an integer $m$ satisfying $\varphi(m) = n$.  In this
section we prove that the problem of deciding whether a set of
integers $\cS$ contains a totient, is {\bf NP}-complete, 
if we assume the following strong form of the {\it Hardy--Littlewood\/}
prime $k$-tuplet conjecture, see~\cite{Bal} for several results in this
direction.

\begin{conjecture}\label{HL_conjecture}
There exists an integer $A > 0$ such that the following
holds:  Suppose that $(M_1 x + a_1)(M_2 x + a_2)$ has
no fixed prime divisors as $x$ runs through the integers, and that
$M_1, M_2 > 0$, and $0 \leq a_i < M_i$ for $i=1$ and $2$.
Then, there exists an $x <\log^A (M_1 M_2 + 1)$
such that both $M_1 x + a_1$ and $M_2 x + a_2$ are prime.
\end{conjecture}

We first note that the decision problem is in {\bf NP}, since if we
let $\cL$ be the language consisting of all finite subsets of 
the natural numbers which contain a totient, then we have:  
For each $\cS \in \cL$, suppose $n \in \cS$ is a totient.
Then, there exists a string $s$, of length $\log^{O(1)} n$, which we can
use to verify that ${\cal S} \in \cL$ in polynomial time, 
namely if $s$ is the prime power factorization of any solution $m$
to $\varphi(m) = n$ 
(and, given the prime power factorization of $m$, it is easy to compute
$\varphi(m)$).  Since we can check whether a number is prime in polynomial
time, and therefore check that $s$ is a legitimate prime power factorization
in time $\log^{O(1)} m$, we conclude that $\cL$ is in {\bf NP}.

The problem which we reduce to our decision problem is the following
variant of the subset sum problem, which   is known to be {\bf NP}-complete.

\begin{quote}
{\sc Partition Problem:} Given   $2k\ge 2$ nonnegative integers
$x_1,\ldots,x_{2k}$, where
     $S = x_1 + \cdots + x_{2k}$ is even,  decide whether
there exist  $1 \le {i_1} < \ldots < {i_k}$
with  $ x_ {i_1} + \cdots +x_{i_k} = S/2$.
\end{quote}

Assuming Conjecture~\ref{HL_conjecture}, we show there is
a polynomial time reduction of the  {\sc Partition}  problem to the
problem of deciding whether there exist integers $m$ satisfying
$\varphi(m) = n$, for a certain small set of values of $n$.

To prove this theorem, we  require the following result which could be of
independent interest.

\begin{theorem}
\label{congruence_theorem}
Given an odd number $k \geq 1$ and given $2k$ integers
$x_1,\ldots,x_{2k}$, we can construct in
polynomial time a series of
congruence classes $a_i \pmod{M}$, $(a_i,M)=1$, such that if
$N_1,\ldots,N_{2k}$ are any numbers satisfying $N_i \equiv a_i \pmod{M}$,
and if $\{i_1,\ldots,i_\ell\} \subset \{1,\ldots,2k\}$, with
$\ell \leq k$, then
\begin{equation} \label{condition0}
\gcd(2N_{i_1} \cdots N_{i_\ell} + 1,\ M)=1 \quad  \iff\  \quad
\ell = k \   \text{and} \  x_{i_1} + \cdots + x_{i_\ell} = S/2;\
\end{equation}
\begin{equation} \label{condition1}
   N_i -1\ \nmid\ 4N_1\cdots N_{2k}, \qquad  i=1,\ldots,2k;
\end{equation}
\begin{equation} \label{condition2}
\gcd(2N_1\cdots N_{2k} + 1,\ M)\ >\ 1, \quad  \text{and}\quad
\gcd(4N_1\cdots N_{2k} + 1,\ M)\ >\ 1.
\end{equation}
\end{theorem}

\begin{proof}  First, we let $R_1,\ldots,R_{k-1}$ be the first
consecutive primes
greater than $k$.  Next, given
$$
A=1 + \sum_{i=1}^{2k} |x_i|,
$$
we let $U_1,\ldots,U_t$ be the first consecutive primes greater than $R_{k-1}$
such that
$$
\prod_{i=1}^t U_i\ >\ 2A.
$$
Finally, we let $v = U_t$, and then let $V_1,\ldots,V_v$ be consecutive primes
greater than $U_t$.  Then, we let
$$
M=8 \cdot 3 \cdot 5 \cdot \prod_{h=1}^{k-1} {2^{R_h} + 1 \over 3}
\prod_{i=1}^t \prod_{j=1}^v {2^{U_i V_j} - 1 \over
(2^{U_i}-1)(2^{V_j} - 1)}.
$$
We claim that this integer $M$ satisfies $\log M = (Ak)^{O(1)}$, which
can be proved by repeated use of the Prime Number Theorem; also,
we claim that each of these factors are coprime to the others, which can be
proved by repeated use of the fact that $(2^H - 1, 2^K - 1) = 2^{(H,K)} - 1$.

We let $a_1,\ldots,a_{2k}$ all be in the same class modulo $8 \cdot 3 \cdot 5$,
defined via the Chinese remainder theorem as follows:
$$
a_i\ \equiv\ 1 \pmod{8};\qquad a_i\ \equiv\ 2 \pmod{3};\qquad a_i\
\equiv\ 4 \pmod{5};
$$
and, for $j=1,\ldots,k-1$, we let
\begin{equation} \label{R_conditions}
a_i\ \equiv\ 2^{g_j} \pmod{(2^{R_j} + 1)/3},
\end{equation}
where for $g_j$ is any solution to $1 + j g_j \equiv R_j \pmod{2R_j}$
(for $j$ odd there is a unique $g_j$; and for $j$ even, there are two
values $g_j$ that satisfy this).

The congruence condition modulo $8$ ensures that~\eqref{condition1} holds;
the congruence modulo $3$ forces the first part of~\eqref{condition2}
to hold; and the condition modulo $5$ forces the second part
of~\eqref{condition2} to hold.
Finally, the condition~\eqref{R_conditions}
   ensures that gcd$(2N_{i_1}\cdots N_{i_\ell} + 1, M) = 1$
implies $\ell = k$, which is part of~\eqref{condition0}.

Now, for $i=1,2,\ldots,t$, we let
$$
\{ \theta(i,1),\ldots,\theta(i,U_i-1)\}=\{ 0,\ldots, U_i - 1\}
\setminus \{ S/2 \pmod{U_i} \};
$$
that is, for every $i=1,2,\ldots,t$, the values of $\theta(i,j)$ run
through the congruence classes
modulo $U_i$, omitting the class $S/2 \pmod{U_i}$.  Next, let
$$
\delta_{i,j}\ \equiv\ k^{-1} \pmod{U_i V_j},\quad  0 \leq \delta_{i,j}
\leq U_i V_j -1.
$$
Then, for $i=1,2,\ldots,t$, $j=1,2,\ldots,U_i-1$, and $\ell =1,\ldots,2k$,
we let
$$
a_\ell\ \equiv\ -2^{V_jx_\ell + \delta_{i,j}(V_j \theta(i,j) - 1)}
\pmod{ {2^{U_i V_j} -1 \over (2^{U_i}-1)(2^{V_j} - 1)}}
$$

Then, if $\{x_{n_1},\ldots,x_{n_k}\}$ is any $k$-element subset of
$k$ of  $\{x_1, \ldots, x_{2k}\}$ such that
$x_{n_1} + \cdots + x_{n_k} \neq S/2$, we must have that for some
$i=1,2,\ldots,t$ and $j=1,2,..,U_i-1$,
$$
x_{n_1} + \cdots + x_{n_k} \equiv \theta(i,j) \pmod{ U_i};
$$
and so, on letting $T = (2^{U_iV_j} - 1)/(2^{U_i} - 1)(2^{V_j} -1)$, we
see that if $N_i \equiv a_i \pmod{M}$, then
\begin{eqnarray*}
2N_{n_1} \cdots N_{n_k} + 1\ &\equiv&\ (-1)^k 2^{1 + V_j(x_{n_1} +
\cdots + x_{n_k} - k \delta(i,j) \theta(i,j)) - k\delta(i,j)} + 1
   \\
&\equiv& -2^{U_i V_j I} + 1 \equiv 0 \pmod{T},
\end{eqnarray*}
where $I$ is some integer.  Conversely, if
$x_{h_1} + \cdots + x_{h_k} = S/2$, then one can show that
$(2N_{h_1} \cdots N_{h_k} + 1, M) = 1$.  Thus, we have
established~\eqref{condition0}, and the result
follows.
\end{proof}

Now are now ready to   prove our main result.

\begin{theorem}\label{np_theorem} Suppose that
$x_1,\ldots,x_{2k}$ is an input of the {\sc Partition} problem.
Let
$$
B = \sum_{i =1}^{2k} \log(x_i +2).
$$
Then,
in polynomial time, we construct a set of $s = B^{O(1)}$ integers
$n_1,\ldots,n_s$ such that the answer to  the corresponding {\sc
Partition} problem
is ``Yes'' if and only if for some $i=1,2,\ldots,s$  we have that
$\varphi(m) = n_i$ has a solution.
\end{theorem}

\begin{proof}
Suppose $x_1,\ldots,x_{2k}$ are given.
We may assume that $k$ is odd, since if $k$ is even, then we can enlarge
our set $\{x_1,\ldots,x_{2k}\}$ by two new elements  $x_{2k+1} = x_{2k+2} = 0$.

Now, suppose that $p_1,\ldots,p_{2k}$ are
a set of primes satisfying $p_i \equiv a_i \pmod{M}$, $p_i > M$.
Then, as a consequence of~\eqref{condition0}, \eqref{condition1},
and~\eqref{condition2} of Theorem~\ref{congruence_theorem},
one can see that if there is a solution $m$ to
$$
\varphi(m)=4p_1\cdots p_{2k},
$$
then $m = P_1P_2$ or $2P_1P_2$, where $P_1$ and $P_2$ are both primes
satisfying
$$
P_1=2p_{i_1} \cdots p_{i_k} + 1,\ \ \ {\rm and\ \ \ } P_2=2p_{j_1}
\cdots p_{j_k} + 1,
$$
where $\{p_{i_1},\ldots,p_{i_k}\} \cup \{p_{j_1},\ldots,p_{j_k}\} =
\{p_1,\ldots,p_{2k}\}$. Moreover, we   have
\begin{equation} \label{x_eqn}
x_{i_1} + \cdots + x_{i_k}=S/2=x_{j_1} + \cdots + x_{j_k}.
\end{equation}

Now suppose that there are two subsets of $\{x_1,\ldots,x_{2k}\}$
satisfying~\eqref{x_eqn}.  Let $\ell$ be one of the numbers $2,3,\ldots,k+2$,
and suppose we are lucky and have
$1 \in \{i_1,\ldots,i_k\}$ and $\ell \in \{j_1,\ldots,j_k\}$,
or have $1 \in \{j_1,\ldots,j_k\}$ and $\ell \in \{i_1,\ldots,i_k\}$;
certainly, for one of these values $\ell=2,3,\ldots,k+2$ this must hold.
We suppose that $1 \in \{i_1,\ldots,i_k\}$ and $\ell \in \{j_1,\ldots,j_k\}$.
Let $\{t_1,\ldots,t_{2k-2}\} = \{1,2,\ldots,2k\} \setminus \{1,\ell\}$.
Then, assuming conjecture~\ref{HL_conjecture} (speicializing to the
case of one linear form, instead of two), we can pick values
$t_1,\ldots,t_{2k-2} < B^{O(1)}$ such that the numbers $a_i + M t_i$
are all prime; moreover, we can pick these  numbers in time $B^{O(1)}$,
by first picking $t_1$, then $t_2$, and so on.

Now, we consider the polynomials
$$
F(x)=2(a_1 + Mx) \prod_{\substack{u \in \{i_1,\ldots,i_k\} \\ u \neq 1}}
(a_u + Mt_u)+ 1,
$$
and
$$
G(y)=2(a_\ell + My) \prod_{\substack{u \in \{j_1,\ldots,j_k\}\\ u \neq \ell}}
(a_u + M t_u)+ 1.
$$
By~\eqref{condition0},  $F(x)$ and $G(y)$ are coprime to $M$ for
all integers $x,y$, and so have no fixed prime divisors;
moreover, $(a_1 + Mx)F(x)$ and $(a_\ell + My)G(y)$ have no fixed prime
divisors.  So, assuming Conjecture~\ref{HL_conjecture}, if we run through
the values $x,y < B^{O(1)}$ that make $a_1 + Mx$ and $a_\ell + My$
both prime, then among these values $x$ and $y$, there must be a choice
which makes $a_1 + Mx, a_\ell + My, F(x),$ and $G(y)$ all prime.
So, we   have a set of primes $p_1,\ldots,p_{2k}$ of the form
$$
p_1=  a_1 + Mx,\qquad  p_\ell=a_\ell + My,
$$
and
$$
p_i  =\ a_i + Mt_i, \qquad i=2,\ldots,\ell-1,\ell+1,\ldots,2k.
$$
These primes satisfy the congruence conditions $p_i \equiv a_i \pmod{M}$.
Furthermore, we also have that $2p_{i_1} \cdots p_{i_k} + 1 =F(x)$ is prime,
as is $2p_{j_1} \cdots p_{j_k} + 1 = G(y)$.  So, if we let
$ n(x,y) = 4p_1\cdots p_{2k}$, then we get a solution
$\varphi(F(x)G(y)) = n(x,y)$.  So, by running through choices for
$x,y < B^{O(1)}$, and $\ell=2,3,\ldots,k+2$, we are guaranteed to
hit upon a value $n(x,y)$ having a solution $\varphi(m) = n(x,y)$,
as long as there is a subset of   $\{x_1,\ldots,x_{2k}\}$  summing to $S/2$.

Conversely, if there is no subset of   $\{x_1,\ldots,x_{2k}\}$
summing to $S/2$, then
either $F(x)$ or $G(y)$ is  an odd composite number, and so  they  fail to
satisfy $\varphi(F(x) G(y)) = n(x,y)$ for all values $x,y$.

Thus, the {\sc Partition}  problem can be reduced, in polynomial time,
to the problem of deciding whether $\varphi(m) = n$ for a set
of $B^{O(1)}$ values $n$,  which finishes the proof.
\end{proof}

\section{Inverting the Euler Function and  Integer Factorisation}
   \label{sec:Factor}

The algorithm of Theorem~\ref{thm:Algorithm} assumes that
the prime number factorisation  of $n$ is given. Here we show
the factorisation problem  for integers from $\cP_2$ can be
reduced in in probabilistic  polynomial
time to the problem of inverting the Euler function.

\begin{theorem}
\label{thm:Factoing}
Given an algorithm that finds
$\Psi(m)$ in time $(\# \Psi^*(m) + \tau(m) + \log m)^{O(1)}$,
without being given the prime
factorisation of $n$,
one can factor  integers $n \in \cP_2$ in probabilistic  polynomial
time.
\end{theorem}

   \begin{proof} Let $\pi(X;r,a)$ denote the number of primes $\ell \le X$ with
$\ell \equiv a \pmod r$.
We need the following result which is a greatly relaxed
version of Theorem~2.1 of~\cite{AGP}.
Namely,  if $r$ is a sufficiently large prime number then for $X\ge r^3$
\begin{equation}
\label{eq:No Sigel}
\pi(X,4r,a) \ge \frac{X}{4r\log X}.
\end{equation}
for any integer $a$ with $\gcd(a,4r) = 1$.

Now, assume we are given sufficiently large  odd $n = pq \in \cP_2$.
We choose two positive  integers $k_1, k_2 \le n^3$ and   consider the product
$4(2k_1+1)(2k_2+1)n$.

It is clear that if  $4(2k_1+1)(2k_2+1)n=\varphi(m)$ then $\omega(m) \le 3$.
More precisely, it is possible only for the values of $m$ of the
form
\begin{enumerate}
\item $m = \ell$ or $m = 2 \ell$ or  $m = 4 \ell$ where $\ell$ is prime;
\item $m = \ell_1\ell_2$  or  $m = 2\ell_1\ell_2$   where $\ell_1,
\ell_2$ are prime;
\end{enumerate}

In each case  of the first group $\ell$ is uniquely defined (and clearly there
are at most two suitable values of $\ell$).

Both cases of the second type occur simultaneously with the same
values of $\ell_1, \ell_2$ which (up to a permutation)  are either of the form
$$
\ell_1 = 2d_1+1 , \qquad \ell_2 =  2d_2pq + 1,
$$
or of the form
$$
\ell_1 = 2d_1p+1 , \qquad \ell_2 = 2 d_2q + 1,
$$
where $d_1$ and $d_2$ are divisors of $(2k_1+1)(2k_2+1)$
with $d_1d_2 = (2k_1+1)(2k_2+1)$.
Therefore, there are at most
$$
2 \tau((2k_1+1)(2k_2+1)) \le 2\tau(2k_1+1)\tau(2k_2+1)
$$
possible solutions of the second kind.
We see from~\eqref{eq:sums of tau} then the total
number of positive integers $k \le X$ with $\tau(k) \ge \log^3 X$
is $O(X\log^{-2} X)$.
Thus from~\eqref{eq:No Sigel} (applied with $r =p$ and $a = 2r +1$)
we derive that there are at
least
$$
\frac{4n^3p}{ 4p \log n^3} + O(n^3\log^{-2} n ) \ge \frac{n^3}{2\log n^3}
$$
   positive  integers $k_1 \le n^3$ for which
simultaneously $2(2k_1+1) p+ 1$ is prime  and $\tau(2k_1+1) \le \log^3 n$.
Similarly, we have at least the same number of positive integers $k_2 \le n^3$
for which simultaneously $2(k_2+1)q+ 1$ is prime  and $\tau(2k_2+1)
\le \log^3 n$.

For each such pair of integers $k_1,k_2$  we
see that the cardinality of $\Psi(4(2k_1+1)(2k_2+1)n)$ is polynomially
bounded,  namely, $\# \Psi(4(2k_1+1)(2k_2+1)n) = O(\log^6 n)$, and
contains a solution
of the form
\begin{equation}
\label{eq:Good Solution}
m = (2(2k_1+1)p+1)(2(2k_2+1)q+1)
\end{equation}
from which, together with the equation $n = pq$, the primes $p$ and $q$
can be trivially found (we certainly have to try all
values of $m \in \Psi(4(2k_1+1)(2k_2+1)n)$ in order
to find the one of the form~\eqref{eq:Good Solution}).

These considerations  naturally lead to the following probabilistic
algorithm which finds the above pair of $k_1, k_2$ and thus
the primes $p$ and $q$.

Assume that the inverting algorithm outputs $\Psi(N)$
in time bounded by $\(\# \Psi(N) \log N\)^A$
for some constant $A>0$.
We choose  integers $k_1,k_2$ uniformly at
random in  the interval $[1, n^3]$ and  use the algorithm
to compute $ \Psi(4(2k_1+1)(2k_2+1)n)$.  If the time
it takes exceeds $\log^{8A} N$  this means that
$\# \Psi(4(2k_1+1)(2k_2+1)n) \ge \log^7 N$ and
we simply terminate the algorithm and choose another
pair $k_1, k_2$. It is clear that in the expected time $O(\log^6 n)$
we find the desired pair of $k_1, k_2$.
\end{proof}


\end{document}